\documentclass[letterpaper,12pt]{article}

\usepackage{color}
\usepackage{amssymb}
\usepackage{amsmath}
\usepackage{amscd}
\usepackage{mathtools}
\usepackage[all]{xy}
\usepackage{tikz-cd}

\usepackage{hyperref}
\hypersetup{colorlinks=true,linkcolor=blue,anchorcolor=blue,citecolor=blue}
\usepackage{cleveref}

\newtheorem{theo}{Theorem}[section]
\newtheorem{prop}[theo]{Proposition}
\newtheorem{lem}[theo]{Lemma}
\newtheorem{cor}[theo]{Corollary}
\newtheorem{defi}[theo]{Definition}
\newtheorem{rem}[theo]{Remark}

\def \deg {{\rm{deg}}}
\def \Br {{\rm{Br}}}

\def \Pic {{\rm {Pic}}}

\def \Gal {{\rm{Gal}}}

\def \A{{\mathbb A}}
\def \P{{\mathbb P}}
\def \Spec {{\rm{Spec}}}

 \def\ov{\overline}

\def \Z {{\mathbb Z}}
\def \Q {{\mathbb Q}}
\def \F {{\mathbb F}}

\def\C{{\mathbb C}}

\def\lra{\longrightarrow}

\def\N{{\rm N}}
\def\H{{\rm H}}

\def\e{\varepsilon}

\def\m{{\bf m}}

\newcommand{\bthe}{\begin{theo}}
\newcommand{\ble}{\begin{lem}}
\newcommand{\bpr}{\begin{prop}}
\newcommand{\bco}{\begin{cor}}
\newcommand{\bde}{\begin{defi}}
\newcommand{\brem}{\begin{rem}}
\newcommand{\ethe}{\end{theo}}
\newcommand{\ele}{\end{lem}}
\newcommand{\epr}{\end{prop}}
\newcommand{\eco}{\end{cor}}
\newcommand{\ede}{\end{defi}}
\newcommand{\erem}{\end{rem}}

\title{Generic diagonal conic bundles revisited}
\author{Alexei N.~Skorobogatov and Efthymios Sofos}

\date{\today}
\begin{document}
\maketitle

\begin{abstract}
We prove a stronger form of our previous result
that Schinzel's Hypothesis holds for
100\% of $n$-tuples of integer polynomials satisfying the usual necessary
conditions, where the primes represented by the polynomials are subject to
additional constraints in terms of Legendre symbols,
as well as upper and lower bounds.
We establish the triviality of 
the Brauer group of generic diagonal conic bundles over the projective line.
Finally, we
give an explicit lower bound for the probability that diagonal conic bundles in 
certain natural families have rational points.
\end{abstract}

\section*{Introduction}

In our previous paper we proved that Schinzel's Hypothesis (H) holds for
100\% of $n$-tuples of integer polynomials satisfying the usual necessary
conditions \cite[Thm.~1.2]{SS}.
Here we give some improvements, complements and further applications of the
results of \cite{SS} relevant to diagonal conic and quadric bundles over the projective line.

In the first section we prove a stronger form of \cite[Thm.~1.2]{SS}
where the primes represented by the polynomials are required to satisfy
additional conditions in terms of Legendre symbols, 
as well as upper and lower bounds,
see Theorem \ref{teor2}.
Using this result, in Corollary \ref{conic} we give a simplified proof of 
a weaker form of the Hasse principle for random diagonal conic bundles over the projective line
\cite[Thm.~6.1]{SS}, with a bound for the least solution. We prove an analogous
statement for diagonal quadric bundles of relative dimension 2, see Corollary \ref{quad}.
(It is well-known that quadric bundles of relative dimension at least 3 over the projective line satisfy
the Hasse principle \cite[Prop.~3.9]{CTSS87}.)

The absence of Brauer--Manin conditions in Corollary \ref{conic} is due to
the triviality of the Brauer group of generic diagonal conic bundles
mentioned in \cite[Remark 6.2]{SS} and proved in the second section of this note, 
see Theorem \ref{main} and Corollary~\ref{coco}.

In the last section we give an explicit lower bound for
the density of pairs of integer polynomials $P_1(t)$, $P_2(t)$ of arbitrary fixed degrees
such that the equation
\begin{equation}
P_1(t)x^2+P_2(t)y^2=z^2 \label{night}
\end{equation}
is soluble in $\Z$. When the degrees of $P_1(t)$ and $P_2(t)$ are large, this density
is close to one third.
We also estimate the height of the smallest integer solution of (\ref{night}).

The authors have been partly supported by the EPSRC New Horizons grant 
``Local-to-global principles for random Diophantine equations" (EP/V019066/1).
We are very grateful to the referees for their thorough reading of the paper and helpful
comments.

\section{Schinzel hypothesis on average with quadratic residue conditions}

Non-constant 
polynomials $P_1(t), \ldots, P_n(t)\in\Z[t]$ are called a {\em Schinzel $n$-tuple} if 
the leading coefficient of each $P_i(t)$ is positive, and for every prime $\ell$
the product $\prod_{i=1}^nP_i(t)$ is not divisible by $t^\ell-t$ modulo $\ell$.

The {\em height} of a polynomial $P(t)\in\Z[t]$ is defined 
as the maximum of the absolute values of the coefficients, and is denoted by $|P|$. The height of
an $n$-tuple of polynomials ${\bf P}=(P_1(t), \ldots, P_n(t))\in(\Z[t])^n$
is defined as $|{\bf P}|=\max_{i=1,\ldots,n}(|P_i|)$

The following result is \cite[Thm.~1.2]{SS} with additional properties (\ref{eq:000}),
(\ref{eq:0}) and (\ref{eq:3}). 
The proof of  (\ref{eq:3}) uses \cite[Prop.~6.5]{SS} based on Heath-Brown's bound for 
character sums \cite[Cor.~4]{HB95}.

 \bthe \label{teor2}
Fix any $(d_1,\ldots,d_n)   \in \mathbb N^n$,  $\varepsilon>0$,  
$M\in \mathbb N$, $m_0\in \Z/M$ 
and ${\bf Q}\in ((\Z/M)[t])^n$ such that 
$ \deg(Q_i)\leq d_i$ and $\gcd(Q_i(m_0), M)=1$ for all $i=1,\ldots, n$. 
For every $(i,j)\in (\mathbb N \cap[1,n])^2$  with $i<j$ let $\epsilon_{ij}\in  \{1,-1\}$. 
Then for $100\%$ of Schinzel $n$-tuples 
${\bf P} \in (\Z[t])^n$ of respective degrees $d_1,\ldots, d_n$ such that
${\bf P}\equiv {\bf Q} \bmod M$,
there exists a natural number $m$ with the following properties:
\begin{align}
&\min\{P_i(m):1\leq i \leq n \}>  |{\bf P}|  (\log |{\bf P}| )^{\varepsilon/2},  \label{eq:000}\\
& m \leq (\log |\mathbf P| )^{n+\varepsilon},  \label{eq:0}\\
&m \equiv m_0 \bmod M,  \label{eq:1}\\
&P_1(m),\ldots,  P_n(m) \text{ are distinct primes}, \label{eq:2}\\
&\text{the Legendre symbol }\left(\frac{P_i(m)}{P_j(m)}\right) \text{ equals }\epsilon_{ij} \ \text{ for all } i<j.
\label{eq:3}
\end{align}  
\ethe
{\em Proof.} Define  $$C_{\bf P}(x):= \sum_{\substack{ m\in \mathbb N \cap[1,x] \\ \eqref{eq:1}-\eqref{eq:3} } } \prod_{s=1}^n \log P_s(m), \quad\quad
\widetilde \theta_{\bf P}(x):= \sum_{\substack{ m\in \mathbb N \cap[1,x] \\ \eqref{eq:1}-\eqref{eq:2} } } \prod_{s=1}^n \log P_s(m).$$
Let us write $\Omega=\{(i,j):1\leq i < j \leq n \}$. For any $ S\subset \Omega$ define 
$$
T_{  S, {\bf P}}(x):= \sum_{\substack{ m\in \mathbb N \cap[1,x] \\ \eqref{eq:1}-\eqref{eq:2} } } 
\prod_{s=1}^n   (\log P_s(m) )  
\prod_{(i,j)\in S}  \left(\frac{P_i(m)}{P_j(m)}\right).$$
We follow the standard convention that a product indexed
by the elements of an empty set is equal to $1$. In particular, we have 
$T_{\emptyset, {\bf P}}(x)=\widetilde \theta_{\bf P}(x)$.
Assuming~\eqref{eq:2},
the following function takes the value $1$ if~\eqref{eq:3} holds and the value $0$ otherwise:
$$ 2^{-\#\Omega}
\prod_{(i,j)\in\Omega} \left( 1+\epsilon_{ij}  \left(\frac{P_i(m)}{P_j(m)}\right)\right)
=
 2^{-\#\Omega}
 \sum_{ S\subset \Omega }
\prod_{(i,j)\in S} \epsilon_{ij}  \left(\frac{P_i(m)}{P_j(m)}\right)
,$$ which leads to 
 \begin{equation}\label{eq:feeding}
C_{\bf P}(x)= 2^{-\#\Omega} \sum_{ S\subset \Omega }
T_{S, {\bf P}}(x)\prod_{(i,j)\in S} \epsilon_{ij} .
\end{equation}
For the rest of the proof we restrict attention only to the range 
 \begin{equation}\label{eq:rangesum}
(\log H)^{A_1 } <x \leq (\log H)^{ A_2},
\end{equation}
where $A_1,A_2$ are arbitrary fixed constants satisfying $n<A_1<A_2$.
By \cite[Eq.~(6.7)]{SS} the term corresponding to $S=\emptyset$ equals 
 \begin{equation}\label{eq:rangesum2}
\frac{\widetilde\theta_{\bf P}(x)}{2^{\#\Omega} }=\frac{\theta_{\bf P}(x)}{2^{\#\Omega} } +O((\log H)^n)
=\frac{\theta_{\bf P}(x)}{2^{\#\Omega} } +O(x^{n/A_1 })
,\end{equation}
where $\theta_{\bf P}(x)$ is defined similarly to $\widetilde\theta_{\bf P}(x)$
by dropping the condition that the primes $P_1(m),\ldots,  P_n(m)$ are necessarily distinct.
Let us define 
$$\texttt{Poly}(H)
:=
\left\{
{\bf P}\in (\Z[t])^n  : \begin{array}{l}
 \deg(P_i)=d_i,
P_i\equiv Q_i\bmod M \text{ for } i=1,\ldots, n, \ |{\bf P}|\leq H
\end{array}
\right\}
.$$ 
We next show that uniformly in the range~\eqref{eq:rangesum}
and for all   $S\neq \emptyset $
one has 
\begin{equation}\label{eq:recvr}
\frac{1}{H^{d+ n}} \sum_{{\bf P} \in \texttt{Poly}(H) } 
  |T_{ S,{\bf P}}(x)| \ll 
x^{    \frac{1}{2} +\frac{n}{2A_1}        }.
\end{equation} 
Fix $h<k$ such that  $(h,k) \in S$.  Letting 
$$\mathcal F_1:= \hspace{-0,2cm} \prod_{\substack{(i,j)\in S  \\  \{i,j\} \cap \{h,k\} =\{h\} } }   \hspace{-0,2cm}  \left(\frac{P_i(m)}{P_j(m)}\right),    \ 
\mathcal F_2:= \hspace{-0,2cm} \prod_{\substack{(i,j)\in S  \\  \{i,j\} \cap \{h,k\} =\{k\} } }  \hspace{-0,2cm}   \left(\frac{P_i(m)}{P_j(m)}\right),   \
\mathcal G:= \hspace{-0,2cm} \prod_{\substack{(i,j)\in S  \\  \{i,j\} \cap \{h,k\} =\emptyset } }   \hspace{-0,2cm}  \left(\frac{P_i(m)}{P_j(m)}\right)
 $$ makes it plain that $T_{S,{\bf P}}$ takes the form of the function $\eta_{\bf P}$ from
\cite[Def.~6.4]{SS}.
This allows us to apply \cite[Prop.~6.5]{SS} to verify~\eqref{eq:recvr}. 
Feeding~\eqref{eq:rangesum2}-\eqref{eq:recvr} into~\eqref{eq:feeding} yields 
$$ \sum_{{\bf P} \in \texttt{Poly}(H) } \frac{ |C_{ \bf P}(x)-2^{-\#\Omega}
 \theta_{\bf P}(x) | }{\#  \texttt{Poly}(H)  } \ll  
x^{ n/A_1 } +x^{\frac{1}{2} +\frac{n}{2A_1}} \ll x^{\frac{1}{2} +\frac{n}{2A_1}} .$$ 
Hence the number of $\bf P $ in $\texttt{Poly}(H)$ with 
$ |C_{ \bf P}(x)-2^{-\#\Omega} \theta_{\bf P}(x)  |>
x^{\frac{1}{2} +\frac{n}{2A_1}} \log x$ is 
$$\leq  \sum_{{\bf P} \in \texttt{Poly}(H) }  \frac{\left |C_{ \bf P}(x)-2^{-\#\Omega}
\theta_{\bf P}(x) \right| } {x^{\frac{1}{2} +\frac{n}{2A_1}} \log x}\ll 
 \frac{\#  \texttt{Poly}(H) }{\log x}=o(\#  \texttt{Poly}(H) ) .$$ 
Therefore, for $100\%$ of $\bf P$ in $\texttt{Poly}(H)$ one has 
 $$ 
|C_{ \bf P}(x)-2^{-\#\Omega} \theta_{\bf P}(x)  | \leq  x^{    \frac{1}{2} +\frac{n}{2A_1}} \log x
,$$ which implies 
$$
C_{\bf P}(x) \geq 
 \frac{2^{-\#\Omega} \beta_0 x}{2(\log \log x )^{d-n}}
-x^{    \frac{1}{2} +\frac{n}{2A_1}        } \log x
>\frac{x}{ \log x},$$
where  we used \cite[Eq.~(4.10)]{SS} for the last deduction. 
Hence, there exists $m\leq x $ satisfying~\eqref{eq:1}-\eqref{eq:3}.
To verify~\eqref{eq:0}, 
we take $ A_2= n+\frac{7 }{10} \varepsilon $ and $A_1= A_2 -\frac{ \varepsilon}{10} $.
Note that $x \leq (\log H)^{A_2}$ 
and that for $100\%$ of ${\bf P}\in  \texttt{Poly}(H)$
one has $H\leq |{\bf P} |^2$. Therefore,   when $|{\bf P}|$ is sufficiently large,
$$m\leq x \leq (\log H)^{A_2}  \ll (\log |{\bf P}|)^{A_2}
 \leq (\log |{\bf P}|)^{n+\varepsilon},$$ 
which proves~\eqref{eq:0}.
To prove~\eqref{eq:000} we note that if the largest integer $m\leq x $ that satisfies \eqref{eq:1}-\eqref{eq:3} is $m_1$
then 
$$
C_{\bf P}(x)= \sum_{\substack{ m\in \mathbb N \cap[1,x] \\ \eqref{eq:1}-\eqref{eq:3} } } \prod_{s=1}^n \log P_s(m)
=
\sum_{m\in \mathbb N \cap[1,m_1] } O((\log H)^n)=O(m_1(\log H)^n)
.$$ For almost all $\mathbf P$ with $|{\bf P}|\leq H$ we have shown that 
$C_{\mathbf P}(x)>\frac{x}{\log x}$. Combined with the upper bound $C_{\bf P}(x) =O(m_1(\log H)^n)$
this gives   
$$m_1\gg \frac{x}{(\log H)^n \log x}\gg \frac{ (\log H)^{A_1-n}}{\log \log H}=\frac{(\log H)^{\frac{3}{5} \varepsilon }}{\log \log H}.$$
Since $P_i(m_1)=\sum_{j=1}^{d_i} c_{ij} m_1^j$ for some integer vector ${\bf c}$ with $|c_{ij }| \leq H$, we see that 
$P_i(m_1)=c_{i{d_i}} m_1^{d_i}+O(H m_1^{d_{i}-1})$. Note that for almost all $\bf P$ with $|{\bf P}| \leq H$ one has 
$\min_{i,j} | c_{ij}| \geq \frac{H}{\log \log H}$, hence, for all $i$ we have 
$$P_i(m_1)= H m_1^{d_{i}-1} \left( \frac{ c_{i{d_i}}}{ H} m_1+O(1) \right) \gg \frac{Hm_1^{d_i}}{\log \log H}
\geq 
\frac{Hm_1 }{\log \log H}
$$ because 
$m_1/\log \log H \geq (\log H)^{A_1-n} (\log \log H)^{-2}\to+\infty$. 
Hence, for all $i$ one has 
$$  P_i(m_1) \gg \frac{Hm_1}{\log \log H}
 \gg H  \frac{(\log H)^{\frac{3}{5}\varepsilon} }{\log \log H }
> |{\bf P}|   \frac{(\log |{\bf P}| )^{\frac{3}{5}\varepsilon} }{\log \log |{\bf P}|  }
.$$ 
In particular, $ \min_{i} P_i(m_1) > |{\bf P}|  (\log |{\bf P}| )^{\varepsilon/2} $, which concludes the proof of~\eqref{eq:000}.
\hfill $\Box$

\subsection*{Generic diagonal conic bundles}

As an application we give a simplified proof of \cite[Thm.~6.1]{SS} with an added value of
a bound for the least solution.
Finite search bounds for Diophantine equations are not well-studied but are nevertheless relevant
to the theory as they are intimately related to Hilbert's 10th problem for $\mathbb Q$.
We give a search bound that is of \textit{polynomial growth} in the size of the coefficients.

We need the following special case of a theorem of Cassels \cite{MR69217}.
\bpr[Cassels]
\label{cas}
If $f_1,f_2,f_3$ are non-zero integers
such that the quadratic form $\sum_{i=1}^3 f_i x_i^2$ represents zero in $\Q$,
then there exists a solution
$(x_1,x_2,x_3)\in \mathbb N^3$ such that 
$$ \max\{ x_i: i=1,2,3\} \leq 40 \max\{|f_i|:i=1,2,3\}.$$ 
\epr

For $\mathbf m=(m_0,m_1,\ldots, m_r)\in\Z^r$ we write $P(t,\mathbf m)$
for the polynomial $\sum_{k=0}^r m_k t^k$.

\bco \label{conic}
Let $n_1, n_2, n_3 $ be integers
such that $n_1>0$, $n_2>0$, and $n_3\geq 0$, and let $n=n_1+n_2+n_3$. Let 
$a_1, a_2, a_3$ be non-zero integers not all of the same sign. 
Let $d_{ij}$ be natural numbers, for $i=1,2,3$ and $j=1,\ldots,n_i$.
Define $d_i=\sum_{j=1}^{n_i}d_{ij}$ and $d=\sum_{i=1}^3 d_i$. 
Let $\mathcal P$ be the set of 
$$\m=(\m_{ij})\in\bigoplus_{i=1,2,3}\bigoplus_{j=1}^{n_i}\Z^{d_{ij}+1} \simeq \Z^{d+n}$$ 
such that the $n$-tuple
$(P_{ij}(t,\m_{ij}))$ is Schinzel.
Let $\mathcal M$ be the set of $\m\in\mathcal P$ such that for each $p|2a_1a_2a_3$ the equation
\begin{equation}
\label{cinque}
a_1 \prod_{j=1}^{n_1}P_{1,j}(t)\,
x_1^2+a_2 \prod_{k=1}^{n_2}P_{2,k}(t)\,x_2^2
+a_3 \prod_{l=1}^{n_3}P_{3,l}(t)\,x_3^2=0,
\end{equation}
has a solution in $\Z_p$ for which the value of each polynomial $P_{ij}(t)$ is a $p$-adic unit.
Then there is a subset $\mathcal M'\subset \mathcal M$ of density $1$ 
such that for every $\m\in\mathcal M'$ the equation $(\ref{cinque})$
  has a solution $(x_1,x_2,x_3,t)\in \mathbb N^4$ with  $$\max\{x_1,x_2,x_3,t\} \leq 
\max\{|a_i|\}  (\log |{\bf{P}}|)^{(n+1)\max\{d_i\}} |{\bf P}|^{\max\{n_i\}}.$$  
The set $\mathcal M'$ has positive density in $\Z^{d+n}$ ordered by height.
\eco
{\em Proof.} 
By absorbing primes into variables $x,y,z$  we can assume that $a_1a_2a_3$ is square-free. 
Write  $M=8a_1a_2a_3$. 
Local solubility of  (\ref{cinque}) with $t=m$
at an odd prime $p|M$ with an additional condition that 
the value of each $P_{ij}(m)$ is a $p$-adic unit depends only on the value of
$P_{ij}(m)$ modulo $p$. For the prime $2$ the same holds modulo $8$. Thus $\mathcal M$
is a finite disjoint union of subsets given by the condition ${\bf P}\equiv{\bf Q}\bmod M$, where
${\bf Q}$ is an $n$-tuple of polynomials in $(\Z/M)[t]$ for which there exists an $m_0\in\Z$
such that (\ref{cinque}) 
with $t=m_0$ has a solution in $\Z_p$
for each $p|M$, and $\gcd(Q_{ij}(m_0),M)=1$. Let us fix such a ${\bf Q}$ and such an $m_0$.
Then for any ${\bf P}\equiv{\bf Q}\bmod M$ and any $m\equiv m_0\bmod M$
the equation (\ref{cinque}) 
 with $t=m$ is soluble in $\Z_p$ for each $p|M$.

Suppose that $p_{ij}:=P_{ij}(m)$, where $i=1,2,3$ and $j=1,\ldots,n_i$, are distinct primes,
where $P_{ij}(t)\equiv Q_{ij}(t)\bmod M$ and $m\equiv m_0\bmod M$.
Thus $p_{ij}\equiv Q_{ij}(m_0)\bmod M$, hence $p_{ij}$ does not divide $M$.
The local solubility of (\ref{cinque}) with $t=m$ at the
primes not dividing $M$ and not equal to one of the
$p_{ij}$ is clear, since the conic has good reduction modulo such a prime.
It remains to show that we can choose $m\equiv m_0\bmod M$
so that (\ref{cinque}) is solvable at each of the primes $p_{ij}=P_{ij}(m)$.

For $i=1,2,3$ define $\pi_i=\prod_{j=1}^{n_i}p_{ij}$. 
Let $\lambda_{ij}\in\F_2$ be such that $$(-1)^{\lambda_{ij}}=\left(\frac{-a_{i'}a_{i''}}{p_{ij}}\right)=
(p_{ij},-a_{i'}a_{i''})_{p_{ij}},$$ where $\{i,i',i''\}=\{1,2,3\}$. Here the middle term is the Legendre symbol and the right hand term is the 
Hilbert symbol. 
By global reciprocity we obtain
\begin{equation}
(-1)^{\lambda_{ij}}=\prod_{p|M} (p_{ij},-a_{i'}a_{i''})_p. \label{ehat} \end{equation}
Define $\widetilde\lambda_{ij}\in\F_2$ as follows:
\begin{itemize}
\item
 $\widetilde\lambda_{1,j}=\lambda_{1j}$, for $j=1,\ldots,n_1$;

\item
$(-1)^{\widetilde\lambda_{2,k}} = (-1)^{\lambda_{2,k}} (\pi_1,p_{2,k})_2$, 
for $k=1,\ldots,n_2$;

\item
$(-1)^{\widetilde\lambda_{3,l}} = (-1)^{\lambda_{3,l}} (\pi_1\pi_2,p_{3,l} )_2$,
for $l=1,\ldots,n_3$.
\end{itemize}
Since $p_{ij}\equiv Q_{ij}(m_0)\bmod M$, we see from (\ref{ehat}) that the $\lambda_{ij}$
depend only on ${\bf Q}$ and $m_0$. Thus the same is true for the $\widetilde\lambda_{ij}$.

\ble \label{lem}
We have $\sum_{i=1,2,3}\sum_{j=1}^{n_i}\widetilde\lambda_{ij}=0$. 
\ele
{\em Proof.}
By assumption, $a_1,a_2,a_3$ are not of the same sign. 
Since $M=8a_1a_2a_3$, from global reciprocity we obtain
$$\prod_{p|M}(-a_1a_2,-a_1a_3)_p=1.$$ 
Local solubility of the conic (\ref{cinque}) at a prime $p$ is equivalent to
$$
(-a_1a_2\pi_1\pi_2,-a_1a_3\pi_1\pi_3)_p=1,$$
so this holds for all $p|M$. Thus the product of 
$$(-a_1a_2,\pi_3)_p(-a_1a_3,\pi_2)_p(-a_2a_3,\pi_1)_p(-1,\pi_1)_p(\pi_1\pi_2,\pi_1\pi_3)_p$$
over all $p|M$ is 1. Using that each $\pi_i$ is coprime to $M$, we see that the contribution
of the last two factors is $(\pi_1,\pi_2)_2(\pi_2,\pi_3)_2(\pi_1,\pi_3)_2$.
From the definition of $\widetilde\lambda_{ij}$ it is immediate that 
this product equals the product of
$(-1)^{\widetilde\lambda_{ij}}$
over all $i$ and $j$. \hfill$\Box$

\medskip

We continue the proof of Corollary \ref{conic}.

Local solubility of the conic (\ref{cinque}) at $p_{ij}$ is equivalent to the condition
\begin{equation}
(-  a_{i'}a_{i''}\pi_{i'}\pi_{i''},p_{ij})_{p_{ij}}=
(-  a_{i'}a_{i''}\prod_{k=1}^{n_{i'}}p_{i'k} \prod_{l=1}^{n_{i''}}p_{i''l},p_{ij})_{p_{ij}}=1
\label{odin}
\end{equation}
where $\{i,i',i''\}=\{1,2,3\}$.  For pairs $(ij)$, $(i'j')$ such that $i<i'$
define $x_{ij,i'j'}\in\F_2$~by
$$(-1)^{x_{ij,i'j'}}=(p_{i'j'},p_{ij})_{p_{ij}}.
$$ 
If $i>i'$, we define $x_{ij,i'j'}$ to be equal to $x_{i'j',ij}$.
Thus we always have $x_{ab,cd}=x_{cd,ab}$.
We observe that (\ref{odin}) is equivalent to the equation
\begin{equation}
\sum_{k=1}^{n_{i'}}x_{ij,i'k}+\sum_{l=1}^{n_{i''}}x_{ij,i''l}=\widetilde\lambda_{ij}.
\label{tri}
\end{equation}
This is clear for $i=1$. For $i=2$ and $i=3$ we prove (\ref{tri}) using 
$$\prod_{k=1}^{n_{i'}}(p_{i'k},p_{ij})_{p_{ij}}=
(\pi_{i'},p_{ij})_2\prod_{k=1}^{n_{i'}}(p_{i'k},p_{ij})_{p_{i'k}},$$
which immediately follows from global reciprocity.

Consider the system of $n=n_1+n_2+n_3$ linear equations (\ref{tri})
in $n_1n_2+n_2n_3+n_1n_3$ variables $x_{ij,i'j'}$.
By assumption we have $n_1>0$ and $n_2>0$, so we have at least one variable, 
namely $x_{11,21}$,
and $n\geq 2$ equations. The sum of the left hand sides
of all equations (\ref{tri}) is zero, and
it is easy to see that the matrix of this linear system has rank $n-1$.
Thus the linear map given by this matrix is surjective onto the subspace
of vectors with zero sum of coordinates in $(\F_2)^n$.
We conclude that the system (\ref{tri}) is solvable 
for arbitrary $\widetilde\lambda_{ij}$ with zero sum. In our case
this holds by Lemma \ref{lem}, so (\ref{tri}) has a solution, say $\e_{ij,i'j'}\in\F_2$,
where $(ij)$ and $(i'j')$ are pairs as above.

Applying Theorem \ref{teor2},
for ${\bf P}\equiv{\bf Q}\bmod M$ in a subset of density 1
we find a natural number $m\equiv m_0\bmod M$,
$m\leq (\log |{\bf P}|)^{n+1/2}$, such that 
the numbers $P_{ij}(m)$ are distinct primes satisfying
$$\left(\frac{P_{i'j'}(m)}{P_{ij}(m)}\right)=(-1)^{\e_{ij,i'j'}},$$
whenever $i<i'$. Then the conic (\ref{cinque})  with $t=m$
is everywhere locally solvable, and so is solvable in $\Z$. 
Furthermore, Proposition \ref{cas}
ensures the existence of a solution $(x_1,x_2,x_3)\in \mathbb N^3$ such that 
$$\max\{x_1,x_2,x_3\} \leq \max\left \{|a_i \mathcal H_i|:i=1,2,3\right\},$$ 
where $\mathcal H_i=40 
\prod_{j=1}^{n_i}P_{ij}(m)$.
Using $m\leq (\log |{\bf P}|)^{n+1/2}$ we obtain 
$$
|\mathcal H_i| \leq 40  \prod_{j=1}^{n_i} \left( (1+d_{ij}) |{\bf P}|m^{d_{ij}} \right) \ll_{d_{ij}} |{\bf P}|^{n_i} (\log |{\bf{P}}|)^{(n+1 /2)d_i},$$
which is at most  $|{\bf P}|^{n_i} (\log |{\bf{P}}|)^{(n+1 )d_i} $ for all sufficiently large $|{\bf P}|$.
\hfill $\Box$

\subsection*{Generic diagonal quadric bundles of relative dimension 2}

Let $a_0,a_1, a_2, a_3$ be non-zero integers not all of the same sign and let
$a=a_0a_1a_2a_3$.
Let $d_1,\ldots,d_n$ be positive integers and let $d=\sum_{i=1}^nd_i$. For $i=1,\ldots,n$ let
$P_i(t,\m_i)=\sum_{j=0}^{d_i}m_{ij}t^j$ be the generic polynomial of degree $d_i$.
Let $S_0, S_1, S_2, S_3$ be subsets of $\{1,\ldots,n\}$. The equation
\begin{equation}
\label{q1}
a_0 \prod_{i\in S_0}P_i(t)\,x_0^2+
a_1 \prod_{i\in S_1}P_i(t)\,x_1^2+a_2 \prod_{i\in S_2}P_i(t)\,x_2^2
+a_3 \prod_{i\in S_3}P_i(t)\,x_3^2=0
\end{equation}
defines a family of quadrics $Q_t$ parametrised by the affine line with coordinate $t$
over the field $\Q(\m)$, where $\m=(\m_1,\ldots,\m_n)$. The generic fibre $Q_\eta$ is a quadric of dimension 2
over $\Q(t,\m)$.
We can multiply (\ref{q1}) by a non-zero element of $\Q(\m)$ and absorb squares
into coordinates $x_i$ without affecting the isomorphism class of $Q_\eta$.
This allows us to assume without loss of generality that each $a_i$ is square-free with no prime
dividing more than two of the $a_i$, and that no element of $\{1,\ldots,n\}$ 
belongs to more than two of the sets $S_i$. Define
$$\delta=a\prod_{i\in S_0}P_i(t)\prod_{i\in S_1}P_i(t)\prod_{i\in S_2}P_i(t)\prod_{i\in S_3}P_i(t).$$

\bco \label{quad}
In the above notation assume that $\delta$ is not a square in $\ov\Q(t,\m)$.
Let $\mathcal P$ be the set of $\m=(\m_i)\in\Z^{d+n}$ such that the $n$-tuple
$(P_i(t,\m_i))$ is Schinzel.
Let $\mathcal M$ be the set of $\m\in\mathcal P$ such that for each $p|2a$ the equation
$(\ref{q1})$
has a solution in $\Z_p$ for which the value of each polynomial $P_i(t)$ is a $p$-adic unit.
Then there is a subset $\mathcal M'\subset \mathcal M$ of density $1$ 
such that for every $\m\in\mathcal M'$ the equation $(\ref{q1})$ has a solution in $\Z$.
The set $\mathcal M'$ has positive density in $\Z^{d+n}$ ordered by height.
\eco
{\em Proof.} We follow the beginning of proof of Corollary \ref{conic}.
Let $M=8a$. Local solubility of (\ref{q1}) at $t=m$
at an odd prime $p|a$ with an additional condition that 
the value of each $P_{ij}(m)$ is a $p$-adic unit depends only on the value of
$P_{ij}(m)$ modulo $p$. For the prime $2$ the same holds modulo $8$. 
Thus $\mathcal M$
is a finite disjoint union of subsets given by the condition ${\bf P}\equiv{\bf Q}\bmod M$, where
${\bf Q}$ is an $n$-tuple of polynomials in $(\Z/M)[t]$ for which there exists an $m_0\in\Z$
such that (\ref{q1}) with $t=m_0$ has a solution in $\Z_p$
for each $p|M$, and $\gcd(Q_{ij}(m_0),M)=1$. Let us fix such a ${\bf Q}$ and such an $m_0$.
Then for any ${\bf P}\equiv{\bf Q}\bmod M$ and any $m\equiv m_0\bmod M$
the equation (\ref{q1}) with $t=m$ is solvable in $\Z_p$ for each $p|M$.

Suppose that $p_i:=P_i(m)$, where $i=1,\ldots,n$, are distinct primes,
where $P_i(t)\equiv Q_i(t)\bmod M$ and $m\equiv m_0\bmod M$.
Thus $p_i\equiv Q_i(m_0)\bmod M$, hence $p_i$ does not divide $M$.
The local solubility of (\ref{q1}) with $t=m$ at the
primes not dividing $M$ and not equal to one of the
$p_i$ is clear, since the quadric has good reduction modulo such a prime.

Suppose that every element of $\{1,\ldots, r\}$ belongs to exactly one of the sets $S_i$, and
every element of $\{r+1,\ldots,n\}$ belongs to two of these sets. 
The condition on $\delta$ implies $r\geq 1$, so the prime $p_1$ belongs to 
exactly one of the sets $S_i$.
The local solubility of (\ref{q1}) with $t=m$ at the primes $p_i$, where $i=1,\ldots,r$, is automatic:
this follows from the fact that the conic obtained by setting $x_j=0$ in (\ref{q1}),
where $i\in S_j$, has good reduction modulo $p_i$.
It remains to show that 
for ${\bf P}\equiv{\bf Q}\bmod M$ in a subset of density 1
we can choose $m\equiv m_0\bmod M$
so that (\ref{q1}) is solvable at each of the primes $p_{r+1},\ldots,p_n$.

To apply Theorem \ref{teor2} we define the values 
$\epsilon_{ij}=\left(\frac{p_i}{p_j}\right)$ for $i<j$, as follows.
The values of $\epsilon_{1,i}=1$ for $i=2,\ldots, r$ are of no importance and can be chosen
arbitrary. For $i=r+1,\ldots, n$ we define 
$$\epsilon_{1,i}=-\left(\frac{a}{p_i}\right)=-(a,p_i)_{p_i}=-\prod_{p|M}(a,p_i)_p.$$
Since $p_i\equiv Q_i(m_0)\bmod M$, we see that this depends only on ${\bf Q}$ and $m_0$.
For $k\geq 2$ we define $\epsilon_{k,l}=1$ for all $l>k$. 
Then (\ref{q1}) with $t=m$ is solvable in $\Z_{p_i}$ for $i=r+1,\ldots, n$, because
up to multiplication by a square the product of all four coefficients is
$ap_1\ldots p_r$, which is a non-square modulo $p_i$.
This implies solubility in $\Z_{p_i}$ by \cite[Ch.~4, Lemma 2.6]{Cas78}.
 
An application of
Theorem \ref{teor2} finishes the proof since the resulting quadric $Q_m$
is everywhere locally solvable, hence solvable in $\Z$. \hfill $\Box$
\medskip

The case when $\delta$ is a square in $\Q(\m)$ can be reduced to the case of 
conic bundles. Indeed, let $C\to\A^1$ be the conic bundle defined by setting $x_0=0$ in (\ref{q1}).
Then every smooth fibre $Q_t$ is isomorphic to $C_t\times C_t$, see 
\cite[Prop.~7.2.4 (c$''$)]{CTS21}, so $Q_t$ has a rational point if and only if $C_t$
has a rational point. Thus Corollary \ref{conic} gives a similar statement for generic diagonal
quadric bundles for which $\delta$ is a square in $\Q(\m)$.

\medskip

When $\delta$ is a square in $\ov\Q(t,\m)$ but not a square in $\Q(\m)$, prime values of
polynomials seem to be insufficient to prove the analogue of Corollary \ref{quad}.
Indeed, consider the following particular case of (\ref{q1}):
\begin{equation}
\label{q2}
P(t)(x_0^2+x_1^2)+x_2^2-2x_3^2=0,
\end{equation}
where $P(t)$ is the generic polynomial of degree $d$.
Since $\delta$ is not a square in $\Q_2$, this equation is solvable in $\Z_2$
for any non-zero 	value of $P(t)$,
see \cite[Ch.~4, Lemma 2.6]{Cas78}. When $P(m)=p$ is prime, equation (\ref{q2})
with $t=m$ is solvable in $\Z_p$ if and only if either
$-1$ or $2$ is a square modulo $p$. Thus solvability in $\Z_p$ does not hold for 100\%
of Schinzel polynomials congruent to $Q(t)$ modulo 8, where $Q(t)$ is an arbitrary polynomial 
such that $Q(m_0)$ is odd for some $m_0$.
(Take $Q(t)$ to be the constant polynomial $3$. Then $p\equiv 3 \bmod 8$, thus neither
$-1$ nor $2$ is a square modulo $p$.)

\section{Brauer group of generic diagonal conic bundles}

In this section we prove
the triviality of the Brauer group of generic diagonal conic bundles as was
mentioned but not proved in \cite[Remark 6.2]{SS}.

Let $n_1,n_2,n_3$ be non-negative integers such that $n_1$ and $n_2$ are positive. 
Suppose that we have positive integers $d_{ij}$, where $i=1,2,3$ and $j=1,\ldots n_i$.
Let $m_{ijk}$ be independent variables, where $i=1,2,3$, $j=1,\ldots n_i$, and 
$k=0,\ldots, d_{ij}$. Write $\m_{ij}=(m_{ijk})$.
Consider the generic polynomials of degree $d_{ij}$:
$$P_{ij}(t,\m_{ij})=\sum_{k=0}^{d_{ij}}m_{ijk}t^k\ \in \ \Z[t,\m_{ij}].$$
Let $K$ be a field of characteristic zero, and let 
$F$ be the purely transcendental extension of $K$ obtained by adjoining all variables
$m_{ijk}$ for $i=1,2,3$, $j=1,\ldots n_i$, and $k=0,\ldots, d_{ij}$.
For $a_1, a_2, a_3\in K^\times$ consider the subvariety $X'\subset\P^2_F\times\A^1_F$
given by (\ref{cinque}). It is easy to see that $X'$ is smooth and geometrically integral,
and the projection to $\A^1_F$ is a proper morphism whose fibres are conics.
There is a natural compactification of $X'\to \A^1_F$ to a smooth projective  surface $X$
with a conic bundle structure $X\to\P^1_K$.  

The main result of this section is 

\bthe\label{main}
Let $n_1>0$, $n_2>0$, $n_3\geq 0$.
Then the natural map $\Br(F)\to\Br(X)$ is an isomorphism.
\ethe

We need to introduce some notation. 
Let $\ov F$ be an algebraic closure of $F$.

Each polynomial $P_{ij}(t)$ is irreducible over $F$, thus
$$F_{ij}:=F[t]/(P_{ij}(t))$$
is a field extension of $F$ of degree $d_{ij}$. We write 
$\N_{ij}\colon F_{ij}\to F$ for the norm map. 
The norm map $\N_{ij}$ gives rise to a homomorphism 
$F_{ij}^\times/(F_{ij}^{\times})^2\to F^\times/(F^{\times})^2$, which we shall denote also
by $\N_{ij}$.

For distinct pairs $(ij)$ and $(rs)$ define $R_{ij,rs}\in F$ as the resultant of 
$P_{ij}(t)$ and $P_{rs}(t)$ considered as polynomials in $F[t]$.
Since $P_{ij}(t)$ and $P_{rs}(t)$ have no common root in $\ov F$,
we have $R_{ij,rs}\neq 0$. We note that
each $R_{ij,rs}$ is an (absolutely) irreducible polynomial in the variables $m_{ijk}$ over $K$,
see \cite[p.~398]{GKZ}.
The polynomials $R_{ij,rs}$ and $R_{i'j',r's'}$ differ by an element of $K^\times$ if and only if
either $i=i'$, $j=j'$, $r=r'$, $s=s'$ or $i=r'$, $j=s'$, $r=i'$, $s=j'$. In the second case we have
$R_{ij,rs}=(-1)^{d_{ij}d_{rs}} R_{rs,ij}$.

To simplify notation, we write
$p_{ij}=m_{i,j,d_{ij}}$ for the leading coefficient of $P_{ij}(t)$.
Write $\N_{ij}(P_{rs})$ for the norm of the image of $P_{rs}(t)$ in $F_{ij}$. Then we have
\begin{equation}
R_{ij,rs}=p_{ij}^{d_{rs}}\N_{ij}(P_{rs})=(-1)^{d_{ij}d_{rs}}p_{rs}^{d_{ij}}\N_{rs}(P_{ij})
\ \in \ F^\times. \label{2}
\end{equation}

For $i=1,2,3$ we define $P_i(t)= \prod_{j=1}^{n_i}P_{ij}(t)$ and
write $p_i=\prod_{j=1}^{n_i}p_{ij}$ for the leading coefficient of $P_i(t)$.
(When $n_3=0$, we write $p_3=1$ for the leading coefficient of the constant polynomial $1$.)
Let $d_i=\sum_{j=1}^{n_i}d_{ij}$ and let $n=n_1+n_2+n_3$.

\medskip

\noindent{\em Proof of Theorem} \ref{main}. 
Let $A\in\Br(F(t))$ be the class of the conic (\ref{cinque}) over $F(t)$.  
Equivalently, $A$ is the class of the quaternion algebra associated to this conic.
Each closed point $M=\Spec(F_M)$ of $\P^1_F$ gives rise to the residue of $A$ 
at this point; this is an element
of $F_M^\times/(F_M^{\times})^2$, see \cite[\S 11.3.1]{CTS21}.
Let $\alpha_{ij}$ be the residue of $A$ at the closed point of $\P^1_F$ given by $P_{ij}(t)=0$. 
By a standard formula \cite[Eq.~(1.18)]{CTS21} 
we see that $\alpha_{ij}$ is the image of
$-a_{i'}a_{i''}P_{i'}(t)P_{i''}(t)$ in $F_{ij}^\times/(F_{ij}^{\times})^2$, where
$\{i,i',i''\}=\{1,2,3\}$. From (\ref{2}) we obtain that 
$\N_{ij}(\alpha_{ij})\in F^\times/(F^\times)^2$ is the class of the following element of $F^\times$:
$$\N_{ij}(-a_{i'}a_{i''}P_{i'}(t)P_{i''}(t))=
(-1)^{d_{ij}}a_{i'}^{d_{ij}}a_{i''}^{d_{ij}}
\prod_{j'=1}^{n_{i'}}p_{ij}^{-d_{i'j'}}R_{ij,i'j'}
\prod_{j''=1}^{n_{i''}}p_{ij}^{-d_{i''j''}}R_{ij,i''j''}.
$$
For any $j'\leq n_{i'}$ this is a product of 
the irreducible polynomial $R_{ij,i'j'}$ and a rational function coprime to $R_{ij,i'j'}$.
In particular, the norm $\N_{ij}(\alpha_{ij})$ is non-trivial,
hence $\alpha_{ij}$ is non-trivial. By assumption $n_1\geq 1$ and $n_2\geq 1$,
thus each $\alpha_{ij}$ is non-trivial.
In particular, $\alpha_{1,1}$ and $\alpha_{2,1}$ are two non-trivial residues of $A$.
It follows that $A$ does not belong to the image of the natural
map $\Br(F)\to \Br(F(t))$, hence  the map
$\Br(F)\to\Br(X)$ is injective by \cite[Lemma~11.3.3]{CTS21}.

The above calculation also shows that all fibres of $X\to\P^1_F$ over the points of $\A^1_F$
are reduced and irreducible. The fibre at infinity is smooth if $d_1,d_2,d_3$ have the same
parity. Let us first assume that this is the case. By \cite[Prop.~11.3.4]{CTS21}, the cokernel of
$\Br(F)\to\Br(X)$ is the homology group of the following complex:
$$
\Z/2\lra(\Z/2)^n\cong\bigoplus_{i=1,2,3}\bigoplus_{j=1}^{n_i}\Z/2\lra F^\times/(F^\times)^2,
$$
where the first map sends the generator $1\in\Z/2$ to $(1,\ldots,1)$, and the second map
sends the generator of the $(i,j)$-summand $\Z/2$ to $\N_{ij}(\alpha_{ij})\in F^\times/(F^\times)^2$. This
sequence is a complex by Faddeev's reciprocity law \cite[Thm.~1.5.2]{CTS21}.

Let $\e_{ij}\in \{0,1\}$, for $i=1,2,3$ and $j=1,\ldots, n_i$, not all of them zero, be such that 
$$\prod_{i=1,2,3}\prod_{j=1}^{n_i}\N_{ij}(\alpha_{ij})^{\e_{ij}}=1\ \in \ F^\times/(F^\times)^2.$$
The factors given by pairs $(ij)$ and $(rs)$ 
contribute the irreducible element $R_{ij,rs}$ to the left hand side and no other factor does.
Thus if $\e_{ij}=1$ for some $i$ and $j$,
then we must have $\e_{rs}=1$ for all $r\neq i$ and all
$s=1,\ldots,n_r$. Repeating the argument, we obtain that $\e_{ij}=1$ 
for all possible values of $i$ and $j$. This proves Theorem \ref{main} 
in the case when $d_1,d_2,d_3$ have the same parity.

Now suppose that $d_1,d_2,d_3$ do not have the same parity. Write $\{1,2,3\}=\{i,i',i''\}$, where
$d_i$ and $d_{i'}$ have the same parity. 
Let $\alpha_\infty$ be the residue of $A$ at infinity. We calculate
that $\alpha_\infty\in F^\times/(F^\times)^2$ is the image of $-a_ia_{i'}p_i p_{i'}\in F^\times$. 
If $i$ or $i'$ is 1, then this is simply divisible by $p_{1,1}$; if $i$ or $i'$ is 2, then this is simply
divisible by $p_{2,1}$. In all cases we conclude that $\alpha_\infty$ is non-trivial.
Thus the fibre of $X\to\P^1_F$ at infinity is singular, reduced and irreducible.
Now \cite[Prop.~11.3.4]{CTS21} says that the map
$\Br(F)\to\Br(X)$ is injective and its cokernel is the homology group of the following complex:
$$
\Z/2\lra(\Z/2)^{n+1}\cong\Z/2\oplus\bigoplus_{i=1,2,3}\bigoplus_{j=1}^{n_i}\Z/2\lra F^\times/(F^\times)^2.
$$
Here the first map sends the generator $1\in\Z/2$ to $(1,\ldots,1)$, and the second map
sends the generator of the $(i,j)$-summand $\Z/2$ to $\N_{ij}(\alpha_{ij})\in F^\times/(F^\times)^2$ and sends the generator of the extra copy of $\Z/2$ (given by the point
at infinity) to $\alpha_\infty$.

Let $\e_{ij}\in \{0,1\}$, for $i=1,2,3$ and $j=1,\ldots, n_i$, and $\e_\infty\in\{0,1\}$,
not all of them zero, be such that 
$$\alpha_\infty^{\e_\infty}\prod_{i=1,2,3}\prod_{j=1}^{n_i}\N_{ij}(\alpha_{ij})^{\e_{ij}}=1\ \in \ F^\times/(F^\times)^2.$$
Since $\alpha_\infty$ is non-trivial, we must have $\e_{ij}=1$ for some pair $(ij)$.
Then the same argument as above implies that $\e_{ij}=1$ for all pairs $(ij)$.
From Faddeev's reciprocity law (or a direct calculation) we know that 
$$\alpha_\infty=\prod_{i=1,2,3}\prod_{j=1}^{n_i}\N_{ij}(\alpha_{ij})\neq 1\ \in \ F^\times/(F^\times)^2,$$
thus we must have $\e_\infty=1$. This finishes the proof of Theorem \ref{main}. \hfill $\Box$

\bco\label{H^1}
Let $n_1>0$, $n_2>0$, $n_3\geq 0$.
Then we have $\H^1(F,\Pic(X_{\ov F}))=0$.
\eco
{\em Proof.}
Since the smooth proper surface $X_{\ov F}$ is rational, we have $\Br(X_{\ov F})=0$
by the birational invariance of the Brauer group \cite[Cor.~6.2.11]{CTS21}. 
Thus we obtain a functorial exact sequence
\begin{equation}
0\to\Pic(X)\to\Pic(X_{\ov F})^{\Gal(\ov F/F)}\to\Br(F)\to \Br(X)\to\H^1(F,\Pic(X_{\ov F})),\label{sara}
\end{equation}
see \cite[Prop.~5.4.2]{CTS21}. If $Y$ is a smooth fibre of $X\to\P^1_F$ above an
$F$-point of $\P^1_F$, then $\H^1(F,\Pic(Y_{\ov F}))=0$, so by \cite[Remark 5.4.3 (2)]{CTS21}
the last map in (\ref{sara}) is surjective.
It remains to apply Theorem \ref{main}. \hfill$\Box$

\medskip

Let $K=\Q$. We can think of 
$\m=(m_{ijk})$ as the coordinates of the affine space $\A^{d+n}_\Q$, where
$n=n_1+n_2+n_3$ and $d=\sum_{ij} d_{ij}$. 
Let $U\subset \A^{d+n}_\Q$ be the Zariski open subset given by the condition that 
the discriminants and the leading coefficients of all the polynomials $P_{ij}(t)$ and
the resultants of all pairs of the polynomials $P_{ij}(t)$ are non-zero.
We have $U\neq\emptyset$, e.g., because $U(\C)\not=\emptyset$.
We note that $X'_\m\subset\P^2_\Q\times\A^1_\Q$ given by (\ref{cinque})
is smooth and geometrically integral when $\m\in U(\Q)$, because no polynomial in $t$
divides more than one coefficient and if it divides some coefficient then it simply divides it.
For $\m\in U(\Q)$ we denote by $X_\m\to\P^1_\Q$ the conic bundle surface over $\Q$
which is a natural smooth compactification of $X'_\m$. Specialisation at $\m\in U(\Q)$
preserves the degrees of the polynomials $P_{ij}(t)$, thus $X_\m$ can be obtained as the
specialisation of the conic bundle surface $X\to\P^1_F$ considered above.
 
\bco \label{coco}
For 100\% of points $\m\in\Z^{d+n}$ ordered by height the natural map $\Br(\Q)\to\Br(X_\m)$ 
is an isomorphism.
\eco
{\em Proof.} 100\% of points $\m\in\Z^{d+n}$ are contained in $U$, and for
such points the surface $X_\m$ and thus the map $\Br(\Q)\to\Br(X_\m)$ are well defined.
For $\m\in U(\Q)$ we have compatible 
specialisation maps 
$${\rm sp}_\m\colon\Pic(X)\to\Pic(X_\m)\quad\text{and}\quad
{\rm sp}_{\m,\ov\Q}\colon\Pic(X_{\ov F})\to\Pic(X_{\m,\ov\Q}),$$ see \cite[\S 3.1, \S 3.4]{H94}. 
Since $X_{\m,\ov\Q}$ is rational, the map ${\rm sp}_{\m,\ov\Q}$ is an isomorphism
by \cite[Prop.~3.4.2]{H94}.
Arguments based on Hilbert's irreducibility theorem 
(see \cite[p.~240]{H94}) imply that there is a hilbertian subset
${\mathcal H}\subset \Q^{d+n}$ such that for $\m\in \mathcal H$ the map 
${\rm sp}_{\m,\ov\Q}$ induces isomorphisms
$$\Pic(X_{\ov F})^{\Gal(\ov F/F)}\tilde\lra\Pic(X_{\m,\ov\Q})^{\Gal(\ov\Q/\Q)}\quad\text{and}\quad
\H^1(F,\Pic(X_{\ov F}))\tilde\lra\H^1(\Q,\Pic(X_{\m,\ov\Q})).$$
Lemma \ref{H^1} implies that for $\m\in{\mathcal H}$ we have $\H^1(\Q,\Pic(X_{\m,\ov\Q}))=0$.
Next, we have a commutative diagram
$$\xymatrix{\Pic(X)\ar@{^{(}->}[r]\ar[d]_{{\rm sp}_\m}&
\Pic(X_{\ov F})^{\Gal(\ov F/F)}\ar[d]^{{\rm sp}_{\m,\ov\Q}}_\cong\\
\Pic(X_\m)\ar@{^{(}->}[r]&\Pic(X_{\m,\ov\Q})^{\Gal(\ov\Q/\Q)}}$$
 The exact sequence (\ref{sara}) and Theorem \ref{main} imply that the top horizontal
map is an isomorphism. It follows that for $\m\in \mathcal H$ the bottom map is also an isomorphism.
This map fits into an exact sequence
$$
0\to\Pic(X_\m)\to\Pic(X_{\m,\ov \Q})^{\Gal(\ov \Q/\Q)}\to\Br(\Q)\to \Br(X_\m)\to
\H^1(\Q,\Pic(X_{\m,\ov \Q})),$$
from which we conclude that for $\m\in \mathcal H$ 
the natural map $\Br(\Q)\to\Br(X_\m)$ is an isomorphism. The complement 
$\Q^{d+n}\setminus{\mathcal H}$ is a thin set \cite[Section 9.2, Prop.~1]{S97}.
The classical theorem of S.D.~Cohen (see, e.g. \cite[Section 13.1, Thm.~1]{S97}) implies that
${\mathcal H}$ contains 100\% of points of $\Z^{d+n}$, when they are
ordered by height. \hfill $\Box$

\section{Explicit densities}

\subsection*{Densities of Schinzel $n$-tuples}

Let $\ell$ be a prime. For $P(t)\in\F_\ell[t]$ we denote by $Z_P(\ell)$ 
the number of zeros of $P(t)$ in $\F_\ell$.
For $n\in\mathbb N$ and ${\bf d}=(d_1,\ldots,d_n) \in \mathbb N^n$, let
$$
\delta_n(\ell,{\bf d}):=\frac{\#\{{\bf P}\in \F_\ell[t]^n: \deg(P_i)\leq d_i, 
Z_{P_1\cdots P_n }(\ell) \neq\ell \}}{\ell^{n+d_1+\ldots+d_n}}
$$
be the density of Schinzel $n$-tuples modulo $\ell$. Write $d=d_1+\ldots+d_n$.
An explicit expression for $\delta_n(\ell,{\bf d})$ is given in \cite[Eq.~(2.6)]{SS}, but
when $\ell$ is large or small compared to the degrees $d_i$
it is easy to calculate $\delta_n(\ell,{\bf d})$ directly. For example, if
$\ell>d$, then $Z_{P_1\cdots P_n }(\ell) \neq\ell$ if and only if each $P_i(t)$
is a non-zero polynomial, thus
$$\delta_n(\ell,{\bf d})=\prod_{i=1}^n\left(1-\frac{1}{\ell^{d_i+1}}\right).$$

\ble \label{tre}
If $ \ell \leq 1+  \min \{d_i: 1\leq i \leq n \}$, then
$$\delta_n(\ell, {\bf d})=1-\left(1-(1-1/\ell)^n \right)^{\ell}.$$ 
\ele
{\em Proof.} The probability for a given polynomial of positive degree to have a non-zero value 
at a given point of $\F_\ell$ is $1-1/\ell$. The product of $n$ polynomials of positive degrees
does not vanish at a given point of $\F_\ell$ with probability $(1-1/\ell)^n$, so it vanishes
with probability $1-(1-1/\ell)^n$. If $\ell\leq d_i$ for all $i$, these events are independent,
so the product of polynomials vanishes everywhere on $\F_\ell$ with probability
$1-\delta_n(\ell, {\bf d})=\big(1-(1-1/\ell)^n \big)^{\ell}$. \hfill $\Box$

\brem \label{rem:222} {\rm
For $\ell=2$ the assumptions of the lemma are
 always met, so we have
$\delta_n(2, {\bf d} ) = 1/2^{n-1}- 1/4^n$. 
In particular, the density of Schinzel $n$-tuples is always at most $1/2^{n-1}$, hence
goes to $0$   when $n\to \infty$. }
\erem

\subsection*{The case $n=2$}
Expression \cite[Eq.~(2.6)]{SS} is complicated 
to evaluate if $1+\min\{d_i\}<\ell\leq d$.
We now give a more practical lower bound in the case $n=2$.

\ble \label{quatro}
We have
$$
\delta_2(\ell, {\bf d} )
\geq 
\left(1-\frac{1}{\ell^{1+d_1}}\right)
\left(1-\frac{1}{\ell^{1+d_2}}\right)
- \frac{2^\ell}{\ell^\ell}.$$
\ele
{\em Proof.}
In the counting function in the definition of $\delta_2(\ell,{\bf d})$ we can assume that neither polynomial
is identically zero, thus 
$$\delta_2(\ell, {\bf d} )=\frac{(\ell^{d_1+1} -1 )(\ell^{d_2+1} -1 )- N}{\ell^{d_1+d_2+2}},$$ where 
$$N=\#\{(P_1,P_2)\in (\F_\ell[t]\setminus \{0\})^2: \deg(P_i)\leq d_i, Z_{P_1P_2}(\ell)=\ell\}.$$
Since $P_1(t)P_2(t)$ vanishes everywhere on $\F_\ell$, we see that 
  $$N\leq \sum_{S\subset \F_\ell}
\#\{(P_1,P_2)\in (\F_\ell[t]\setminus \{0\})^2: \deg(P_i)\leq d_i, \ P_1(S)=P_2(S^c)=0  \},$$ 
where $S^c=\F_\ell\setminus S$. 
Since $\#S \leq d_1 $ and $\#S^c\leq d_2$, we obtain
$$\begin{aligned}N\leq&
\sum_{\substack{ S\subset \F_\ell \\ \ell-d_2  \leq \#S \leq d_1 }}
\#\{(P_1,P_2)\in (\F_\ell[t]\setminus \{0\})^2: \deg(P_i)\leq d_i, \ P_1(S)=P_2(S^c)=0\}\\
= &
\sum_{\substack{ S\subset \F_\ell \\ \ell-d_2  \leq \#S \leq d_1 }}
 (-1+\ell^{1+d_1-\#S})
 (-1+\ell^{1+d_2-\ell+\#S}) 
\leq \ell^{2+d_1+d_2-\ell}\,2^\ell,\end{aligned}$$
which gives the desired bound. \hfill $\Box$

\subsection*{Equation $P_1(t)x^2+P_2(t)y^2=z^2$}

The main goal of this section is to give an explicit lower bound for
the density of pairs of integer polynomials $P_1(t)$, $P_2(t)$ of fixed degrees $d_1$, $d_2$,
respectively, such that the above equation is solvable in $\Z$. 

We use Corollary \ref{conic} with $n_1=n_2=1$, $n_3=0$, $a_1=a_2=1$, and $a_3=-1$.
We only need to deal with the prime $\ell=2$. 
If $a$ and $b$ are odd integers, then the Hilbert symbol
$(a,b)_2$ equals 1 if and only if $a$ or $b$ is $1\bmod 4$. 
Thus we need to calculate
the probability $\sigma_2$ that a pair of polynomials $P_1(t), P_2(t)\in(\Z/4)[t]$
of degrees $\deg(P_1)\leq d_1$, $\deg(P_2)\leq d_2$ satisfies the following conditions:
\begin{itemize}
\item[(a)]
$P_1(0)$ and $P_2(0)$ are both odd, or $P_1(1)$ and $P_2(1)$ are both odd, and 
\item[(b)]
$P_i(n)\equiv 1\bmod 4$ for some $i\in\{1,2\}$ and some $n\in\Z/4$.
\end{itemize}
Condition (a) is the Schinzel condition at $2$, and condition (b) is
the triviality of the Hilbert symbol at $2$. Recall that $d=d_1+d_2$. We have 
$$
\sigma_2=\delta_2(2, {\bf d} )-4^{-d-2}\#T=1/2 - 1/4^2-4^{-d-2}\#T,
$$
where $T$ is the set of pairs $P_1(t), P_2(t)\in(\Z/4)[t]$
of degrees $\deg(P_1)\leq d_1$, $\deg(P_2)\leq d_2$ satisfying condition (a) and
taking only values $0$, $2$, and $3$ modulo $4$.

There are $2^{e-1}$ polynomials of degree at most $e$ in $\F_2[t]$ with given values 
at $0$ and $1$. Each of these can be lifted to exactly $2^{e-3}$ polynomials of degree at most $e$
in $(\Z/4)[t]$ with given values at the elements of $\Z/4$ that are compatible with 
the values at $0$ and $1$ modulo $2$. Thus there are $4^{e-2}$ polynomials $P(t)\in(\Z/4)[t]$
of degree at most $e$ with given values at the elements of $\Z/4$ such that the values at 
0 and 2 have the same parity, and similarly for the values at 1 and 3.

For $i\in\Z/4$ we write $v_i=P_1(i)$
and $v'_i=P_2(i)$. Condition (a) is equivalent to the condition that $v_0v'_0$ is odd
or $v_1v'_1$ is odd.
We obtain that $\#T$ is the product of $4^{d-4}$ and the number of 
$(v_j)  \in \{0,2,3\}^4$, $(v'_j)  \in \{0,2,3\}^4$ such that
$$v_0v'_0 \text{ or } v_1v'_1\, \text{is odd, and}\,
v_0-v_2, v_1-v_3, v'_0-v'_2, v_1'-v_3'\ \text{are even}.$$ 
Since $2\mid v_0-v_2$, the values of $(v_0,v_2)$ can only be 
\begin{equation}\label{eq:lll}
(0,0), (0,2), (2,0), (2,2), (3,3).
\end{equation}
The same can be said about $(v_1,v_3), (v'_0,v'_2) , (v'_1,v'_3) $. 
If $v_0v'_0$ is odd, then $v_0=v_2=v'_0=v'_2=3$.
Then there are 25 possibilities for $(v_1,v_3, v'_1,v'_3)$ since both $ (v_1,v_3)$ and $(v'_1,v'_3)$ can be chosen arbitrarily from the list in~\eqref{eq:lll}.
A similar counting gives another $25$ cases when $v_1v'_1$ is odd.
We have counted twice the cases where both  $v_0v'_0$ and $v_1v'_1$ are odd.
This implies that all $v_i$ and $v'_i$ equal 3, so there is in fact only one such case.
We obtain $\#T=25+25-1=49$, and hence
$$\sigma_2=1/2 - 1/4^2-\frac{49}{4^{6}} =\frac{1743}{4096}=42.553\ldots.$$

Evaluating the product $\sigma_2\prod_{\ell\geq 3}\delta_2(\ell, {\bf d})$ using Lemmas
\ref{tre} and \ref{quatro} gives the following bound.

\bpr
Let $d_1$ and $d_2$ be positive integers. The density of pairs of integer polynomials $P_1(t)$, $P_2(t)$ of  degrees $\deg(P_1)=d_1$ and $\deg(P_2)=d_2$ with positive leading
coefficients, ordered by height,
such that the equation
$$P_1(t)x^2+P_2(t)y^2=z^2 $$
has a solution $(x,y,z,t)\in \mathbb{N}^4$ with $\max\{x,y,z,t\}\leq (\log |\mathbf P|)^{3\max\{d_i\}}|\mathbf P|$, is at least 
$$\frac{1743}{4096}
\prod_{3 \leq \ell\leq 1+\min d_i  } 
\left(1-\left(  \frac{2\ell-1}{\ell^2}\right)^{\ell}\right)
\prod_{ \ell >  1+\min d_i  }  \left(\left(1-\frac{1}{\ell^{1+d_1}}\right)
\left(1-\frac{1}{\ell^{1+d_2}}\right)
+O\left( \frac{2^\ell}{\ell^\ell}\right)\right).$$ 
\epr

When $\min\{d_i\} \to \infty $ this product becomes asymptotic to 
$$\mathcal P:=\frac{1743}{4096}\prod_{\text{prime } \ell \geq 3  } 
\left(1-\left(2/\ell-1/\ell^2\right)^{\ell}\right).$$
This is strictly exceeding $0.35$. Verifying $\mathcal P>0.35$ with any 
computer package   will  lead to troubles related to accumulating 
rounding errors owing to     $(2/\ell-1/\ell^2)^{\ell}$. Thus, it is better to use the bound 
$$ \left(2/\ell-1/\ell^2\right)^{\ell} <(2/\ell)^\ell \leq 1/\ell^4, \quad \text{for} \ \ \ell \geq 17,$$  which implies that  
$$\mathcal P_0 \prod_{\text{prime } \ell \geq 17   } \left(1-1/\ell^4 \right)<\mathcal P <
\mathcal P_0 ,$$ where  $$\mathcal P_0 =\frac{1743}{4096} \prod_{\text{prime } \ell \in [3,13]  } \left(1-\left(2/\ell-1/\ell^2\right)^{\ell}\right)
.$$  We can easily calculate that  $$\prod_{ \ell \geq 17  }  (1-1/\ell^4 ) = 0.9999723\ldots
\text{ and }  \mathcal P_0 =0.3504\ldots ,$$ therefore, we have $0.3503< \mathcal P<0.3504$.

Department of Mathematics, 
South Kensington Campus,
Imperial College London,
SW7~2AZ United Kingdom \  \ and \  \ 
Institute for the Information Transmission Problems,
Russian Academy of Sciences,
Moscow, 127994 Russia

\texttt{a.skorobogatov@imperial.ac.uk}

\bigskip

Department of Mathematics,
University of Glasgow,
University Place,
Glasgow,  G12~8QQ United Kingdom

\texttt{efthymios.sofos@glasgow.ac.uk}

\end{document}